\documentclass[12pt]{article}
\usepackage{amssymb}

\newtheorem{definition}{Definition}[section]
\newtheorem{theorem}[definition]{Theorem}
\newtheorem{lemma}[definition]{Lemma}
\newtheorem{corollary}[definition]{Corollary}
\newtheorem{remark}[definition]{Remark}

\newtheorem{note}[definition]{Note}

\def\qed{{~~~\vrule height .75em width .4em depth .2em}}

\def\Z{\mathbb Z}

\newcommand{\beast}{\begin{eqnarray*}}
\newcommand{\eeast}{\end{eqnarray*}}
\newcommand{\spn}{\mathop{\rm span}\nolimits}
\newcommand{\diag}{\mathop{\rm diag}\nolimits}

\begin{document}

\title{ \bf  Three Mutually Adjacent Leonard Pairs}

\author{Brian Hartwig {\footnote{
Department of Mathematics, University of
Wisconsin, 480 Lincoln Drive, Madison WI 53706-1388 USA}
}}

\date{}
\maketitle
\begin{abstract}
Let $\Bbb K$ denote a field of characteristic $0$ and let $V$ denote a
vector space over $\Bbb K$ with positive finite dimension.  Consider
an ordered pair of linear transformations $A:V \rightarrow V$ and
$A^*:V \rightarrow V$ that satisfies both conditions below:
\begin{enumerate}
\item
There exists a basis for $V$ with respect to which the matrix
representing $A$ is diagonal and the matrix representing $A^*$ is
irreducible tridiagonal.
\item
There exists a basis for $V$ with respect to which the matrix
representing $A^*$ is diagonal and the matrix representing $A$ is
irreducible tridiagonal.
\end{enumerate}
We call such a pair a {\it Leonard pair} on $V$.  Let $(A,A^*)$ denote
a Leonard pair on $V$.
A basis for $V$ is said to be {\it standard} for $(A,A^*)$
whenever it satisfies $(i)$ or $(ii)$ above.  A basis for $V$ is said
to be {\it split} for $(A,A^*)$ whenever with respect to
this basis the matrix representing
one of $A,A^*$ is lower bidiagonal and the matrix representing the
other is upper bidiagonal.  Let $(A,A^*)$ and $(B,B^*)$ denote
Leonard pairs on $V$.  We say these pairs are {\it adjacent} whenever
each basis for $V$ which is standard for
$(A,A^*)$ (resp. $(B,B^*)$) is split for $(B,B^*)$ (resp. $(A,A^*)$).
Our main results are as follows.

\medskip
\noindent
{\bf Theorem 1}  {\it There exist at most $3$ mutually adjacent
Leonard pairs on $V$ provided the dimension of $V$ is at least $2$.}

\medskip
\noindent
{\bf Theorem 2}  {\it Let $(A,A^*)$, $(B,B^*)$, and $(C,C^*)$ denote
three mutually adjacent Leonard pairs on $V$.  Then for each of these
pairs, the eigenvalue sequence and dual
eigenvalue sequence are in arithmetic progression.}
\medskip

\noindent
{\bf Theorem 3}  {\it Let $(A,A^*)$ denote a Leonard pair on $V$ whose
eigenvalue sequence
and dual eigenvalue sequence are in arithmetic progression.  Then
there exist Leonard pairs $(B,B^*)$
and $(C,C^*)$ on $V$ such that $(A,A^*)$, $(B,B^*)$, and $(C,C^*)$ are
mutually adjacent.}
\medskip

\end{abstract}

\medskip


\section{Leonard Pairs}

Throughout the paper, $\Bbb K$ will denote a field of characteristic
$0$ and $V$ will denote a vector space over $\Bbb K$ with positive
finite dimension. 

We begin by
recalling the notion of a Leonard pair
\cite{ter1},\cite{ter2},\cite{ter3},\cite{ter4},\cite{ter5},\cite{ter6},\cite{ter7},\cite{ter8}.
We will use
the following terms.  Let $M$ denote a square matrix.  Then $M$ is called
$tridiagonal$ whenever each nonzero entry lies on either the
diagonal, the subdiagonal, or the superdiagonal.  Assume $M$ is
tridiagonal.  Then $M$ is called $irreducible$ whenever each entry on
the subdiagonal is nonzero and each entry on the superdiagonal is nonzero.

\begin{definition}  
\label{def:1st}
\cite{ter1}
\rm
By a $Leonard$ $pair$ on $V$, we mean an
  ordered pair $(A,A^*)$, where $A : V \rightarrow V$ and $A^* : V
  \rightarrow V$ are linear transformations that satisfy both $(i)$
  and $(ii)$ below:
\begin{enumerate}
\item
There exists a basis for $V$ with respect to which the matrix
representing $A$ is diagonal and the matrix representing $A^*$ is
irreducible tridiagonal.
\item
There exists a basis for $V$ with respect to which the matrix
representing $A^*$ is diagonal and the matrix representing $A$ is
irreducible tridiagonal.
\end{enumerate}
\end{definition}

\begin{note}
\rm
It is a common notational
convention to use $A^*$ to represent the conjugate-transpose of $A$.
We are {\it not} using this convention.  In a Leonard pair $(A,A^*)$,
the linear transformations $A$ and $A^*$ are arbitrary subject to
$(i)$ and $(ii)$ above.  
\end{note}

In this paper we introduce the notion of adjacency for Leonard pairs.
Our main results are summarized as follows.  We show that there exist
at most three mutually adjacent Leonard
pairs on $V$ provided that $\dim V \ge 2$.  Given
three mutually adjacent Leonard pairs on $V$, we show that the
eigenvalue sequence
and dual eigenvalue sequence of each is in arithmetic progression.
Given a Leonard pair on $V$ whose eigenvalue
sequence and dual eigenvalue sequence are in arithmetic progression,
we show that there exist two additional Leonard pairs on $V$ such that
all three Leonard pairs are mutually adjacent.

For the rest of this section we recall some basic results concerning
Leonard pairs.

\begin{definition}
\rm
Let $(A,A^*)$ denote a Leonard pair on $V$.
Let $W$ denote a vector space over $\Bbb K$ with positive finite
dimension and let $(B,B^*)$ denote a
Leonard pair on $W$.  By an {\it isomorphism of Leonard pairs} from
$(A,A^*)$ to $(B,B^*)$, we mean an isomorphism of vector spaces
$\sigma:V \rightarrow W$ such that $\sigma A \sigma^{-1} = B$ and
$\sigma A^* \sigma^{-1} = B^*$.  We say $(A,A^*)$ and $(B,B^*)$ are
{\it isomorphic} whenever there exists an isomorphism of Leonard pairs from
$(A,A^*)$ to $(B,B^*)$.
\end{definition} 

\begin{lemma}
\label{th:aaes}
\cite[Lemma 1.3]{ter1}
Let $(A,A^*)$ denote a Leonard pair on $V$.  Then the eigenvalues of $A$
(resp. $A^*$) are mutually distinct and contained in $\Bbb K$.  
\end{lemma}

\begin{lemma}
\cite[Lemma 3.3]{ter1}
\label{th:irred}
Let $(A,A^*)$ denote a Leonard pair on $V$.  Then there does not exist
a proper nonzero subspace $W$ of $V$ such that $AW \subseteq W$ and
$A^*W \subseteq W$. 
\end{lemma}

By a {\it decomposition} of $V$ we mean a sequence $V_0,V_1,\ldots,
V_d$ of one dimensional subspaces of $V$ such that
$$V=V_0+V_1+ \cdots +V_d \qquad (\mbox{direct sum}). $$
Let $v_0,v_1,\ldots, v_d$ denote a basis for $V$ and let
$V_0,V_1, \ldots, V_d$ denote a decomposition of $V$.  We say
$V_0,V_1, \ldots, V_d$ is $induced$ by
$v_0,v_1,\ldots, v_d$ whenever $V_i = \spn (v_i)$ for $0 \leq i \leq d$.

\begin{definition}  
\rm
Let $(A,A^*)$ denote a Leonard pair on $V$.  A basis for $V$ is said
to be {\it $A$-standard} (resp. {\it $A^*$-standard})
whenever with respect to this basis
the matrix representing $A$ (resp. $A^*$) is diagonal and the matrix
representing $A^*$ (resp. $A$) is irreducible tridiagonal.  A
decomposition of $V$ is said to be {\it $A$-standard} (resp. {\it
$A^*$-standard}) whenever it is induced by an $A$-standard basis
(resp. $A^*$-standard basis).  
A basis (resp. decomposition) for $V$ is said to be {\it standard}
for $(A,A^*)$ whenever it is either $A$-standard or $A^*$-standard.  
\end{definition}

Let $a_0, a_1, \ldots, a_n$ be a finite sequence.  By the {\it
inversion} of $a_0, a_1, \ldots, a_n$ we mean the sequence $a_n,
a_{n-1}, \ldots, a_0$.  

Let $(A,A^*)$ denote a Leonard pair on $V$ and let $V_0, V_1,
\ldots, V_d$ denote a decomposition of $V$.  Observe that the
inversion
$V_d, V_{d-1}, \ldots, V_0$ is also a decomposition of $V$.  One
easily verifies that $V_0, V_1, \ldots, V_d$ is $A$-standard
(resp. $A^*$-standard) if and only if the inversion $V_d, V_{d-1},
\ldots, V_0$ is $A$-standard (resp. $A^*$-standard).
Moreover, by \cite[p. 388]{ter 9} there is no
other $A$-standard (resp. $A^*$-standard) decomposition of $V$. 

\section{Flags}

In this section we will discuss the notion of a {\it standard flag} for a
Leonard pair. 

\begin{definition}  
\rm
By a {\it flag} on $V$, we mean a sequence $F_0,F_1, \ldots, F_d$
of subspaces of $V$ such that $F_{i-1} \subset F_i$ for
$1 \leq i \leq d$, $F_i$ has dimension $i+1$ for $0 \leq i \leq d$,
and $F_d = V$.  We call $F_i$ the $i^{\mathrm{th}}$ {\it component} of
the flag.
\end{definition}

The following construction yields a flag on $V$.  Let $V_0, V_1,
\ldots, V_d$ denote a decomposition of $V$.  Set
$$F_i = V_0 + V_1 + \cdots + V_i$$
for $0 \leq i \leq d$.  Observe that the sequence
$F_0,F_1, \ldots, F_d$ is a
flag on $V$.  We say $F_0,F_1,\ldots, F_d$ is {\it induced} by $V_0,
V_1, \ldots, V_d$.

For each Leonard pair, we will define a set of flags as follows.

\begin{definition}
\label{def:lsf}
\rm
Let $(A,A^*)$ denote a Leonard pair on $V$.  A flag on $V$ is
said to be {\it $A$-standard} (resp. {\it $A^*$-standard}) whenever it
is induced by an $A$-standard (resp. $A^*$-standard) decomposition of
$V$.  A flag on $V$
is said to be {\it standard} for $(A,A^*)$ whenever it is either
$A$-standard or $A^*$-standard.  We define $\mathcal F (A,A^*)$ to be the
set of the flags on $V$ which are standard for $(A,A^*)$.
\end{definition}

\begin{lemma}
\label{lem:distn}
Let $(A,A^*)$ denote a Leonard pair on $V$.  No flag on $V$ is both
$A$-standard and $A^*$-standard provided $\dim V \geq 2$.
\end{lemma}

{\it Proof.}  Assume $\dim V \geq 2$.  
Suppose there exists a flag on $V$ that is both $A$-standard and
$A^*$-standard.  For this flag the $0^{\mathrm{th}}$ component is a
one dimensional
subspace of $V$ that is invariant for $A$ and $A^*$.  This contradicts
Lemma \ref{th:irred}.
\qed

\begin{corollary}
\label{cor:ki}
Let $(A,A^*)$ denote a Leonard pair on $V$.  Then $(i)$--$(iii)$ below
are true. 
\begin{enumerate}
\item
The number of $A$-standard flags on $V$ is two if $\dim V \geq 2$ and one if
$\dim V =1$.
\item
The number of $A^*$-standard flags on $V$ is two if $\dim V \geq 2$ and one if
$\dim V =1$.
\item
$| \mathcal F (A,A^*) | = 4$ if $\dim V \geq 2$ and  $| \mathcal F (A,A^*) | = 1$
if $\dim V =1$.
\end{enumerate}
\end{corollary}

{\it Proof.}  Assume $\dim V \geq 2$; otherwise the result is trivial.

Recall from the last paragraph of section 2 that there are exactly two
$A$-standard decompositions of
$V$ and these decompositions are inversions of each other.
The two $A$-standard flags induced by these decompositions are
distinct since their $0^{\mathrm{th}}$ components are distinct.
A similar argument shows there are two $A^*$-standard flags.  By
Lemma \ref{lem:distn} no flag on $V$ is both $A$-standard and
$A^*$-standard, so $| \mathcal F (A,A^*) | = 4$.
\qed

\medskip

We now discuss the notion of {\it opposite} flags.  Let $F_0,F_1,
\ldots, F_d$ and $G_0,G_1, \ldots, G_d$ denote flags on $V$.
These flags are said to be $opposite$ whenever 
$$ F_i \cap G_j = 0 \hbox{ if } i+j < d, \quad 0 \leq i,j \leq d.$$

The following construction produces an ordered pair of opposite flags
on $V$.  Let $V_0, V_1, \ldots, V_d$ denote a decomposition of $V$.
Set $$F_i = V_0 + V_1 + \cdots + V_i,$$
$$G_i = V_d + V_{d-1} + \cdots + V_{d-i}$$
for $ 0 \leq i \leq d$.  Observe that the sequences $F_0,F_1, \ldots, F_d$ and
$G_0,G_1,\ldots, G_d$ are opposite flags on $V$.

Given an ordered pair of opposite flags on $V$, the following construction
produces a decomposition of $V$.  Let $F_0,F_1, \ldots, F_d$ and $G_0,G_1,
\ldots, G_d$ denote an ordered pair of opposite flags on $V$.  Set
$$V_i = F_i \cap G_{d-i}$$
for $0 \leq i \leq d$.  One easily verifies that $V_0, V_1, \ldots,
V_d$ is a decomposition of $V$.

\begin{remark}
\label{rm:of}
\rm
Let $D$ denote the set of all decompositions of $V$, and
let $F$ denote the set of all ordered pairs of opposite
flags on $V$.  In the previous two paragraphs, we defined a map from
$D$ to $F$ and a map from $F$ to $D$.  It is routine to show that these
maps are inverses of one another.  In particular, each of these maps
is a bijection.
\end{remark}

We will use the following notation.

\begin{definition}
\label{def:makedec}
\rm
Let $f,g$ denote an ordered pair of opposite flags on $V$.  Set 
$$V_i = F_i \cap G_{d-i}, \quad 0 \leq i \leq d$$  
where $F_j$ (resp. $G_j$) denotes the $j^{\mathrm{th}}$ component of $f$
(resp. $g$) for $0 \leq j \leq d$.  Since $f$ and
$g$ are opposite, $V_0, V_1, \ldots, V_d$ is a decomposition of $V$.
We denote this decomposition by $[fg]$.
\end{definition}

We now return our attention to Leonard pairs.

\begin{theorem}
\cite[Theorem 7.3]{ter2}
\label{th:mo}
Let $(A,A^*)$ denote a Leonard pair on $V$.  Then the
flags in $\mathcal F (A,A^*)$ are mutually opposite.
\end{theorem}

We will find the following result useful.

\begin{corollary}
\label{cor:di}
Let $(A,A^*)$ denote a Leonard pair on $V$.  Let $x,y$ denote distinct
flags on $V$.  Then the following are equivalent.
\begin{enumerate}
\item 
Each of $x,y$ is $A$-standard.
\item
The flags $x,y$ are opposite and $[xy]$ is an $A$-standard
decomposition.
\end{enumerate}
\end{corollary}

{\it Proof.}  $(i) \Rightarrow (ii)$:  Observe that $x,y$ are distinct
elements of $\mathcal F (A,A^*)$, so by Theorem \ref{th:mo}, $x,y$ are
opposite.  Since $x$ is an $A$-standard flag, by Definition
\ref{def:lsf} there exists an $A$-standard decomposition $V_0,V_1,
\ldots, V_d$ that induces $x$.  Similarly there exists an $A$-standard
decomposition that induces $y$.  This decomposition must be $V_d,V_{d-1},
\ldots, V_0$ by Corollary \ref{cor:ki}$(i)$ and since $x \not= y$.
Observe that the
decomposition $[xy]$ is equal to $V_0,V_1, \ldots, V_d$ and is
therefore $A$-standard.

$(ii) \Rightarrow (i)$:  Combine Remark \ref{rm:of} and Definition
\ref{def:lsf}. 
\qed

\section{The Split Decomposition}

In this section we discuss the split decompositions for a Leonard
pair.  We will use the following terms.  Let $M$ denote a square matrix.
We say $M$ is {\it lower bidiagonal} whenever each nonzero entry
lies on either the diagonal or the subdiagonal.  We say $M$ is {\it
upper bidiagonal} whenever the transpose of $M$ is lower bidiagonal. 

\begin{definition}
\rm
Let $(A,A^*)$ denote a Leonard pair on $V$.  A basis for $V$ is said
to be {\it $LU$-split} for
$(A,A^*)$ whenever with respect to this basis the matrix representing $A$ is
lower bidiagonal and the matrix representing $A^*$ is upper
bidiagonal.  A decomposition of $V$ is said to be {\it $LU$-split} for
$(A,A^*)$ whenever it is induced by a basis for $V$ that is $LU$-split
for $(A,A^*)$.
\end{definition}

\begin{lemma}
\cite[Theorem 4.6]{ITT}
\cite[Corollary 7.6]{ter3}
\label{lem:asd}
Let $(A,A^*)$ denote a Leonard pair on $V$.  Let $V_0,V_1, \ldots,
V_d$ denote a decomposition of $V$.  Then the following are equivalent.
\begin{enumerate}
\item
$V_0,V_1, \ldots, V_d$ is $LU$-split for $(A,A^*)$.
\item
There exist an $A^*$-standard flag $x$ and an $A$-standard flag $y$
such that $[ xy ]$ is equal to $V_0,V_1,\ldots,V_d$. 
\end{enumerate}
\end{lemma}

\begin{definition}
\rm
Let $(A,A^*)$ denote a Leonard pair on $V$.  A basis for $V$ is said
to be {\it $UL$-split} for
$(A,A^*)$ whenever with respect to this basis the matrix representing $A$ is
upper bidiagonal and the matrix representing $A^*$ is lower
bidiagonal.  A decomposition of $V$ is said to be {\it $UL$-split} for
$(A,A^*)$ whenever it is induced by a basis for $V$ that is $UL$-split
for $(A,A^*)$.
\end{definition}

Let $(A,A^*)$ denote a Leonard pair on $V$.  One easily verifies that a
decomposition (resp. basis) of $V$ is $LU$-split for $(A,A^*)$ if and
only if the inversion of that decomposition (resp. basis) is $UL$-split
for $(A,A^*)$.  By this and Lemma \ref{lem:asd} we
obtain the following result.

\begin{lemma}
\label{lem:asdb}
Let $(A,A^*)$ denote a Leonard pair on
  $V$.  Let $V_0,V_1, \ldots, 
V_d$ denote a decomposition of $V$.  Then the following are equivalent.
\begin{enumerate} 
\item
$V_0,V_1, \ldots, V_d$ is $UL$-split for $(A,A^*)$.
\item
There exist an $A$-standard flag $x$ and an $A^*$-standard flag $y$ such that $[ xy ]$
is equal to $V_0,V_1, \ldots,V_d$.
\end{enumerate}
\end{lemma}

\begin{definition}
\rm
Let $(A,A^*)$ denote a Leonard pair on $V$.  A basis
(resp. decomposition) of $V$ is said to be $split$ for $(A,A^*)$
whenever it is either $LU$-split or $UL$-split for $(A,A^*)$. 
\end{definition}

\section{Adjacent Leonard Pairs}

In this section we define what it means for two Leonard pairs on $V$
to be {\it adjacent}.  We begin with a lemma.

\begin{lemma}
\label{th:doalp}
Let $(A,A^*)$ and $(B,B^*)$ denote Leonard pairs on $V$.  Then the
following are equivalent.
\begin{enumerate}
\item
Each decomposition of $V$ that is standard for $(A,A^*)$ is split for
$(B,B^*)$.
\item
Each decomposition of $V$ that is standard for $(B,B^*)$ is split for
$(A,A^*)$. 
\end{enumerate}
\end{lemma}

{\it Proof.}  Assume $\dim V \geq 2$; otherwise the result is trivial.

$(i) \Rightarrow (ii)$:  Consider an $A$-standard decomposition of $V$.
By assumption this decomposition is split for $(B,B^*)$.  
Inverting this decomposition if necessary, we can assume that it is
$UL$-split for $(B,B^*)$.  By Lemma \ref{lem:asdb}, there exists a
$B$-standard flag $x$ and a
$B^*$-standard flag $y$ such that $[xy]$ is this $A$-standard decomposition.
By Corollary \ref{cor:di}, we find $x,y$ are the
$A$-standard flags.
Consider an $A^*$-standard decomposition of $V$.  By assumption this
decomposition is split for $(B,B^*)$. 
Inverting this decomposition if necessary, we can assume that it is
$UL$-split for $(B,B^*)$.  By Lemma \ref{lem:asdb}, there exists a
$B$-standard flag $w$ and a
$B^*$-standard flag $z$ such that $[wz]$ is this $A^*$-standard
decomposition.  By Corollary \ref{cor:di}, we find $w,z$ are the
$A^*$-standard flags.
By Lemma \ref{lem:distn} and since $\dim V \geq 2$, no flag on $V$ is
both $A$-standard and $A^*$-standard.  Therefore $w,x,y,z$ are distinct.
We now see that $x,w$ are the $B$-standard flags and $y,z$ are the
$B^*$-standard flags.  Apparently the decompositions of $V$ that are
standard for $(B,B^*)$ are $[ xw ],[ wx ],[ yz ],[ zy ]$.  By Lemma 
\ref{lem:asd} and \ref{lem:asdb}, each of these is split for $(A,A^*)$.

$(ii) \Rightarrow (i)$: Reverse the roles of $(A,A^*)$ and $(B,B^*)$ in
the proof of $(i) \Rightarrow (ii)$.
\qed

\medskip

Rephrasing Lemma \ref{th:doalp} in terms of bases, we obtain the
following result.

\begin{corollary}
\label{cor:doalp}
Let $(A,A^*)$ and $(B,B^*)$ denote Leonard pairs on $V$.  Then the
following are equivalent.
\begin{enumerate}
\item
Each basis for $V$ that is standard for $(A,A^*)$ is split for
$(B,B^*)$.
\item
Each basis for $V$ that is standard for $(B,B^*)$ is split for
$(A,A^*)$. 
\end{enumerate}
\end{corollary}

\begin{definition}
\rm
\label{def:ad1}
Let $(A,A^*)$ and $(B,B^*)$ denote
Leonard pairs on $V$.  We say these pairs are {\it adjacent} whenever
they satisfy conditions $(i)$ and $(ii)$ in Lemma
\ref{th:doalp} (equivalently, they satisfy
conditions $(i)$ and $(ii)$ in Corollary \ref{cor:doalp}.)
\end{definition}

\begin{lemma}
\label{lem:ca}
Let $(A,A^*)$ and $(B,B^*)$ denote adjacent Leonard pairs on $V$.
Then each of the following $(i)$--$(iv)$ are adjacent Leonard pairs on
$V$.
\begin{enumerate}
\item 
$(A,A^*)$ and $(B,B^*)$.
\item
$(A^*,A)$ and $(B,B^*)$.
\item
$(A,A^*)$ and $(B^*,B)$.
\item
$(A^*,A)$ and $(B^*,B)$.
\end{enumerate}
\end{lemma}

{\it Proof.}  Observe that a basis for $V$ is standard (resp. split)
for $(A,A^*)$ if and only if that basis is standard (resp. split) for
$(A^*,A)$.  A similar statement applies for $(B,B^*)$ and $(B^*,B)$.
The result follows. 
\qed

\medskip

Our next goal is to show that there exist at
most three mutually adjacent Leonard pairs on $V$ provided $\dim V
\geq 2$.  To do this, we first introduce some notation.  

Let $(A,A^*)$ denote a Leonard pair on $V$ and assume $\dim V \geq 2$.
We define a relation $\sim$ on the set
$\mathcal F (A,A^*)$ as follows.  Let $x,y \in \mathcal F (A,A^*)$.  Then $x
\sim y$ whenever
either $x$ and $y$ are both $A$-standard flags or $x$ and $y$ are both
$A^*$-standard flags.  We observe that $\sim$ is an equivalence relation on
$\mathcal F (A,A^*)$.  The relation $\sim$ partitions $\mathcal F (A,A^*)$ into
two equivalence classes, each containing two elements.  We call
$\sim$ the {\it principal relation} induced by $(A,A^*)$.   
 
\begin{lemma}
\label{lem:rty}
Assume $\dim V \geq 2$.  Let $(A,A^*)$
and $(B,B^*)$ denote Leonard pairs on $V$.  Then the
following are equivalent.
\begin{enumerate}
\item
$(A,A^*)$ and $(B,B^*)$ are adjacent.
\item
$\mathcal F (A,A^*) = \mathcal F (B,B^*)$ and the principal relation
induced by $(A,A^*)$ is different from the principal relation induced
by $(B,B^*)$.
\end{enumerate}
\end{lemma}

{\it Proof.}  
$(i) \Rightarrow (ii)$:  By Definition \ref{def:ad1} we find Lemma
\ref{th:doalp}$(i)$ holds.  We argue as in the proof of $(i) \Rightarrow
(ii)$ from Lemma \ref{th:doalp}.  Using the same notation as in that
proof, we see that $\mathcal F (A,A^*)$ and $\mathcal F (B,B^*)$ are both equal
to $\{ w,x,y,z \}$.  Recall that $x,w$ are the $B$-standard flags, that
$y,z$ are the $B^*$-standard flags, that $x,y$ are the $A$-standard
flags, and that $w,z$ are the $A^*$-standard flags.  Therefore the
principal relation induced by $(A,A^*)$ is different from the
principal relation induced by $(B,B^*)$.

$(ii) \Rightarrow (i)$:  The decompositions of $V$ that are standard
for $(A,A^*)$ are split for $(B,B^*)$.  Therefore $(A,A^*)$ and
$(B,B^*)$ are adjacent by Definition \ref{def:ad1}.
\qed

\medskip

We now present our first main result.

\begin{theorem}
Assume $\dim V \geq 2$.  Then there exist at
most three mutually adjacent Leonard pairs on $V$.
\end{theorem}

{\it Proof.}  There are three ways to partition a four element set
into two sets, each of size two.  The result follows from this and Lemma
\ref{lem:rty}.
\qed

\section{The Eigenvalue And Dual Eigenvalue Sequences For A Leonard Pair}

In this section we discuss the eigenvalues of a Leonard pair.

\begin{definition}
\rm
Let $(A,A^*)$ denote a Leonard pair on $V$ and let $V_0,V_1, \ldots,
V_d$ denote an $A$-standard decomposition of $V$.  Recall that for $0
\leq i \leq d$, $V_i$ is an eigenspace for $A$; let $\theta_i$ denote the
corresponding eigenvalue.  We call $\theta_0, \theta_1, \ldots,
\theta_d$ the {\it eigenvalue sequence} for $(A,A^*)$ that corresponds to
$V_0,V_1, \ldots, V_d$.  
\end{definition}

\begin{definition}
\rm
Let $(A,A^*)$ denote a Leonard pair on $V$ and let $V_0^*,V_1^*, \ldots,
V_d^*$ denote an $A^*$-standard decomposition of $V$.  For $0 \leq i \leq
d$, recall $V_i^*$ is an eigenspace for $A^*$; let $\theta_i^*$ denote the
corresponding eigenvalue.  We call $\theta_0^*, \theta_1^*, \ldots,
\theta_d^*$ the {\it dual eigenvalue sequence} for $(A,A^*)$ that
corresponds to $V_0^*,V_1^*, \ldots, V_d^*$.
\end{definition}

Let $(A,A^*)$ denote a Leonard pair on $V$.  Observe that if
$\theta_0, \theta_1, \ldots, \theta_d$ is an eigenvalue sequence for
$(A,A^*)$ then so is $\theta_d, \theta_{d-1}, \ldots,
\theta_0$ and there is no other eigenvalue sequence for $(A,A^*)$.  A
similar result holds for the dual eigenvalue sequences of $(A,A^*)$.

We recall a basic property of the eigenvalue and dual eigenvalue
sequences.  

\begin{lemma}
\cite[Theorem 11.1]{ITT}
\label{th:estuff}
Let $(A,A^*)$ denote a Leonard pair on $V$.
Let $\theta_0, \theta_1, \ldots , \theta_d$ (resp. $\theta_0^*,
\theta_1^*, \ldots , \theta_d^*$) denote an eigenvalue
sequence (resp. dual eigenvalue sequence) for $(A,A^*)$. 
Then the scalars
$$ \frac{\theta_{i-2}-\theta_{i+1}}{\theta_{i-1}-\theta_{i}}, \qquad 
 \frac{\theta_{i-2}^*-\theta_{i+1}^*}{\theta_{i-1}^*- \theta_{i}^*}$$
are equal and independent of $i$ for $2\leq i\leq d-1$.  
\end{lemma}

Parametric expressions for the eigenvalue sequences and dual eigenvalue
sequences of a Leonard pair can be found in \cite[Theorem 11.2]{ITT}.

Let $(A,A^*)$ denote a Leonard pair on $V$.
Let $\theta_0, \theta_1, \ldots , \theta_d$ (resp. $\theta_0^*,
\theta_1^*, \ldots , \theta_d^*$) denote an eigenvalue
sequence (resp. dual eigenvalue sequence) for $(A,A^*)$.  In this
paper we will encounter the special case in which the scalars  
$$\frac{\theta_{i}-\theta_{i+1}}{\theta_{i-1}
  -\theta_{i}},  \qquad \frac{\theta_{i}^*-\theta_{i+1}^*}{\theta_{i-1}^*
  -\theta_{i}^*},$$
are equal and independent of $i$ for $1 \leq i \leq d-1$.  The next
three lemmas prepare us for this special case.

\begin{lemma}  
\label{lem:seq2}
Let $d$ denote a nonnegative integer and let $\theta_0,\theta_1,
  \ldots, \theta_d$ denote a sequence of
  mutually distinct scalars in $\Bbb K$.  Given $q \in \Bbb K$ such
  that $q \not= 0$ and $q \not= 1$, the following are
  equivalent.
\begin{enumerate}
\item  For $1 \leq i \leq d-1$,
$$ \frac{\theta_{i}-\theta_{i+1}}{\theta_{i-1}
  -\theta_{i}} =q.$$
\item
There exists $\alpha,\beta \in \Bbb K$ such that $\alpha \not= 0$ and
  $\theta_i = \alpha q^i +\beta$ for $0 \leq i \leq d$.
\end{enumerate}
\end{lemma}

{\it Proof.}  $(i) \Rightarrow (ii)$:  Observe $\theta_{i+1} -
(1+q) \theta_{i} + q \theta_{i-1} =0$ for $1 \leq i \leq d-1$.  The
characteristic polynomial of this recursion is $x^2 -
(1+q) x + q = 0$ and this polynomial has roots at $x=1$ and $x=q$.  We
conclude that there exists $\alpha,\beta \in \Bbb K$ such that
$\theta_i = \alpha q^i + \beta$ for $0 \leq i \leq d$.  Furthermore,
$\alpha \not= 0$ since $\theta_0,\theta_1, \ldots,  \theta_d$ are
mutually distinct.

$(ii) \Rightarrow (i)$:  This direction is clear.
\qed

\begin{definition}
\label{def:has}
\rm
Let $d$ denote a nonnegative integer and let $\theta_0,\theta_1,
\ldots, \theta_d$ denote a sequence of
mutually distinct scalars in $\Bbb K$.  Given $q \in \Bbb K$ such
  that $q \not= 0$ and $q \not= 1$, we call this sequence
{\it q-classical} whenever it satisfies the equivalent conditions
$(i)$, $(ii)$ from Lemma \ref{lem:seq2}.
\end{definition}

\begin{note}
\label{n:ote}
\rm
Referring to Definition \ref{def:has}, assume the sequence $\theta_0,
\theta_1, \ldots, \theta_d$ is $q$-classical.  Then $q^i \not= 1$ for
$1 \leq i \leq d$.  
\end{note}

{\it Proof.}  Immediate from Lemma \ref{lem:seq2}$(ii)$ and the fact
that $\theta_0, \theta_1, \ldots, \theta_d$ are mutually distinct.
\qed

\begin{note}
\rm
Referring to Definition \ref{def:has}, the sequence $\theta_0,
\theta_1, \ldots, \theta_d$ is $q$-classical if and only if the sequence
$\theta_d,\theta_{d-1}, \ldots, \theta_0$ is $q^{-1}$-classical.
\end{note}

\begin{lemma}  
\label{lem:seq1}
Let $d$ denote a nonnegative integer and let $\theta_0,\theta_1,
  \ldots,  \theta_d$ denote a sequence of
  mutually distinct scalars in $\Bbb K$.  The following are equivalent.
\begin{enumerate}
\item  For $1 \leq i \leq d-1$,
$$ \frac{\theta_{i}-\theta_{i+1}}{\theta_{i-1}
  -\theta_{i}} =1.$$
\item
There exists $\alpha,\beta \in \Bbb K$ such that $\alpha \not= 0$ and
  $\theta_i = \alpha i +\beta$ for $0 \leq i \leq d$.
\end{enumerate}
\end{lemma}

{\it Proof.}  Routine.
\qed

\begin{definition}
\label{def:had}
\rm
Let $d$ denote a nonnegative integer and let $\theta_0,\theta_1,
\ldots, \theta_d$ denote a sequence of
mutually distinct scalars in $\Bbb K$.  We call this sequence
{\it arithmetic} whenever it satisfies the equivalent conditions
$(i)$, $(ii)$ from Lemma \ref{lem:seq1}.
\end{definition}

\begin{note}
\rm
Referring to Definition \ref{def:had}, the sequence $\theta_0,
\theta_1, \ldots, \theta_d$ is arithmetic if and only if the sequence
$\theta_d,\theta_{d-1}, \ldots, \theta_0$ is arithmetic.
\end{note}

We now return our attention to Leonard pairs.

\begin{lemma} 
\label{lem:eval}
Let $(A,A^*)$ denote a Leonard pair on $V$.  Let $\theta_0,\theta_1,
\ldots,  \theta_d$ denote an eigenvalue sequence for $(A,A^*)$.
Assume $d \geq 3$ and
\begin{eqnarray}
\frac{\theta_{1}-\theta_{2}}{\theta_0
  -\theta_{1}} = \frac{\theta_{2}-\theta_{3}}{\theta_1
  -\theta_{2}}.  \label{eq:terw}
\end{eqnarray}
Then 
$$\frac{\theta_{i}-\theta_{i+1}}{\theta_{i-1}
  -\theta_{i}}$$ is independent
of $i$ for $1 \leq i \leq d-1$.
\end{lemma}

{\it Proof.}  We show
\begin{eqnarray}
\frac{\theta_{1}-\theta_{2}}{\theta_0
  -\theta_{1}} = \frac{\theta_{i}-\theta_{i+1}}{\theta_{i-1} 
  -\theta_{i}}  \label{eq:ind}
\end{eqnarray}
for $ 1 \leq i \leq d-1$.  We proceed by induction.  Let $i$ be given.
If $i=1$ then (\ref{eq:ind}) holds so assume $2 \leq i \leq d-1$.  
By Lemma \ref{th:estuff}, 
$$\frac{\theta_0-\theta_3}{\theta_1-\theta_2} =
\frac{\theta_{i-2}-\theta_{i+1}}{\theta_{i-1}
  -\theta_{i}}.$$  
Using this and induction we find
$$- \frac{\theta_{0}-\theta_{1}}{\theta_1 -\theta_{2}} + 
\frac{\theta_{0}-\theta_{3}}{\theta_1 -\theta_{2}} - 1
=
- \frac{\theta_{i-2}-\theta_{i-1}}{\theta_{i-1} -\theta_{i}} + 
\frac{\theta_{i-2}-\theta_{i+1}}{\theta_{i-1}-\theta_{i}}
- 1.$$
Reducing we get 
$$\frac{\theta_{2}-\theta_{3}}{\theta_1-\theta_{2}} = 
\frac{\theta_{i}-\theta_{i+1}}{\theta_{i-1}-\theta_{i}}.$$
Evaluating this equation using (\ref{eq:terw}) we obtain
(\ref{eq:ind}).  The result follows.
\qed

\section{Eigenvalue And Dual Eigenvalue Sequences For Adjacent Leonard Pairs}

In this section we will discuss the eigenvalues and dual eigenvalues
of adjacent Leonard pairs.  We use the following notation.

\begin{definition}
\label{def:s-u}
\rm
Let $(A,A^*)$ and $(B,B^*)$ denote adjacent
Leonard pairs on $V$, and assume $\dim V \geq 2$.  
Recall by Lemma \ref{lem:rty} that $\mathcal F (A,A^*) = \mathcal F (B,B^*)$
and the principal relations induced by $(A,A^*)$ and $(B,B^*)$ are
distinct.  Let $w,x,y,z$ denote the elements of $\mathcal F (A,A^*) =
\mathcal F (B,B^*)$, ordered so that the flag types are given as follows.
$$
\begin{tabular}{c|cc}
 & $B$-standard flags & $B^*$-standard flags \\
\hline
$A$-standard flags & $w$ & $x$ \\
$A^*$-standard flags & $z$ & $y$  \\
\end{tabular}
$$
With respect to this labeling,
\begin{enumerate}
\item
Let $\theta_0, \theta_1, \ldots , \theta_d$ denote the eigenvalue
sequence for $(A,A^*)$ associated with the decomposition $[wx]$.
\item
Let $\theta_0^*, \theta_1^*, \ldots , \theta_d^*$ denote the dual eigenvalue
sequence for $(A,A^*)$ associated with the decomposition $[yz]$.
\item
Let $\eta_0, \eta_1, \ldots , \eta_d$ denote the eigenvalue
sequence for $(B,B^*)$ associated with the decomposition $[zw]$.
\item
Let $\eta_0^*, \eta_1^*, \ldots , \eta_d^*$ denote the dual eigenvalue
sequence for $(B,B^*)$ associated with the decomposition $[xy]$.
\end{enumerate}
\end{definition}

\begin{lemma} 
\label{lem:th}
With reference to Definition \ref{def:s-u},
\begin{eqnarray}
 \frac{(\theta_{d-i} - \theta_d)(\theta_{d-i}-\theta_{d-1}) \cdots
  (\theta_{d-i}-\theta_{d-j+1})}{(\theta_{d-j}-\theta_d)
  (\theta_{d-j}-\theta_{d-1}) \cdots
  (\theta_{d-j}-\theta_{d-j+1})} =
\frac{(\eta_0-\eta_{j+1})(\eta_0-\eta_{j+2}) \cdots 
(\eta_0-\eta_{i})}{(\eta_{j}-\eta_{j+1})
  (\eta_{j}-\eta_{j+2})  
  \cdots (\eta_{j}-\eta_{i})} \label{eq:ttot}
\end{eqnarray}
for $0 \leq j \leq i \leq d$.
\end{lemma}

{\it Proof.}  Let $V_0,V_1, \ldots, V_d$ denote the decomposition
$[wx]$ and let $U_0,U_1, \ldots, U_d$ denote the decomposition $[zw]$.
Let $0 \not= u \in U_0$.  Observe that $u$ is an eigenvector for
$A^*$ with eigenvalue
$\theta_d^*$.  Define $u_i = (A-\theta_d) \cdots (A-\theta_{d-i+1}) u$
for $0 \leq i \leq d$.  By \cite[p. 841]{ter2}, $u_i$ is a
basis for $U_i$ for $0 \leq i \leq d$.  Moreover, $u_0, u_1, \ldots,
u_d$ is a basis for V.  By \cite[Section 19]{ter3}, there exists a basis
$v_0,v_1, \ldots, v_d$ for $V$ such that $v_i \in V_i$ for $0 \leq
i \leq d$ and $\displaystyle{\Sigma^d_{i=0}}v_i = u$.  Let $T_1$
denote the transition matrix from the basis $v_d, v_{d-1}, \ldots , v_0$
to the basis $u_0, u_1, \ldots, u_d$.  By \cite[Theorem 15.2]{ter2}, $T_1$
is lower triangular with entries
$$T_1(i,j) =  (\theta_{d-i}-\theta_d) \cdots (\theta_{d-i}-\theta_{d-j+1})
\qquad 0 \leq j \leq i \leq d.$$

Observe that $v_d$ is an eigenvector for $B^*$ with eigenvalue
$\eta_0^*$.  Define  $v_i' =
(B-\eta_0) \cdots (B-\eta_{d-i-1})v_d$ for $0 \leq i \leq d$.
By \cite[p. 841]{ter2}, $v_i'$ is a basis for $V_i$ for $0
\leq i \leq d$.  Moreover, $v_0', v_1', \ldots, v_d'$ is a basis for V.
By \cite[Section 19]{ter3}, there exists a basis $u_0',u_1',
\ldots, u_d'$ for $V$ such that $u_i' \in U_i$ for $0 \leq
i \leq d$ and $\displaystyle{\Sigma^d_{i=0}}u_i' = v_d$.
Let $T_2$ denote the transition matrix from the basis
$v_d', v_{d-1}', \ldots , v_0'$ to the basis $u_0', u_1', \ldots, u_d'$.  By
\cite[Theorem 15.2]{ter2}, $T_2$ is lower triangular with entries
$$T_2(i,j) = \frac{1}{(\eta_j-\eta_0) \cdots (\eta_j-\eta_{j-1})}
\frac{1}{(\eta_j-\eta_{j+1}) \cdots (\eta_j-\eta_i)} \qquad 0 \leq j
\leq i \leq d.$$ 

Let $D_1$ denote the transition matrix from the basis $v_d, v_{d-1},
\ldots, v_0$ to the basis $v_d', v_{d-1}', \ldots , v_0'$.  Since $v_i, v_i'
\in V_i$ for $0 \leq i \leq d$, we find $D_1$ is diagonal.   For $0
\leq i \leq d$, let $\alpha_i$ denote the $(i,i)$ entry of $D_1$.
Observe that $v_d' = v_d$, so $\alpha_0=1$.  Let $D_2$ denote the
transition matrix from the basis $u_0', u_1', \ldots, u_d'$ to the basis $u_0,
u_1, \ldots , u_d$. Since $u_i, u_i' \in U_i$ for $0 \leq i \leq d$,
we find $D_2$ is diagonal.  For $0 \leq i \leq d$, let $\beta_i$
denote the $(i,i)$ entry of $D_2$. 

Observe that $D_1 T_2 D_2$ is the transition matrix from the basis $v_d,
v_{d-1}, \ldots, v_0$ to the basis $u_0, u_1, \ldots, u_d$.  This
transition matrix is also given by $T_1$, so $T_1=D_1 T_2 D_2$.
Pick $i,j \in \Z$  with $0 \leq j \leq i \leq d$.  Equating the
$(i,j)$ entries of $T_1$ and $D_1 T_2 D_2$ we find
\begin{eqnarray}
(\theta_{d-i}-\theta_d) \cdots (\theta_{d-i}-\theta_{d-j+1})
= \alpha_i \beta_j \frac{1}{(\eta_j-\eta_0) \cdots (\eta_j-\eta_{j-1})}
\frac{1}{(\eta_j-\eta_{j+1}) \cdots (\eta_j-\eta_i)}. 
\label{eq:dtd}
\end{eqnarray}
Setting $i=0$, $j=0$ in (\ref{eq:dtd}), we find that $\alpha_0 \beta_0
= 1$, and therefore $\beta_0 =1$.  Setting $j=0$ in
(\ref{eq:dtd}), we find
\begin{eqnarray}
\alpha_i = (\eta_0-\eta_1) \cdots
(\eta_0-\eta_i) \qquad \qquad 0 \leq i \leq d 
\label{eq:ali}.
\end{eqnarray}   
Setting $j=i$ in (\ref{eq:dtd}) and using (\ref{eq:ali}), we find 
\begin{eqnarray}
\beta_i = (\theta_{d-i}-\theta_d) \ldots
(\theta_{d-i}-\theta_{d-i+1}) \frac{(\eta_i-\eta_0) \cdots
  (\eta_i-\eta_{i-1})}{(\eta_0-\eta_1) \cdots (\eta_0-\eta_i)} \qquad
\qquad 0 \leq i \leq d 
\label{eq:bei}.
\end{eqnarray}
Evaluating (\ref{eq:dtd}) using (\ref{eq:ali}) and (\ref{eq:bei}), we
get (\ref{eq:ttot}) for $0 \leq j \leq i \leq d$.
\qed

\begin{lemma}  
\label{lem:vis3}
With reference to Definition \ref{def:s-u}, the scalars
\begin{eqnarray}
\frac{\theta_{i} -\theta_{i+1}}{\theta_{i-1}-\theta_{i}},
\qquad \qquad
\frac{\eta_{i}-\eta_{i+1}}{\eta_{i-1} -\eta_{i}}
\label{eq:star}
\end{eqnarray}
are equal and independent of $i$ for $1 \leq i \leq d-1$.
\end{lemma}

{\it Proof.}  Assume $d \geq 2$; otherwise the result is trivial.

For $1 \leq i \leq d-1$, let $\epsilon_i$ (resp. $\delta_i$) denote
the fraction on the left (resp. right) in (\ref{eq:star}).  We show
$\epsilon_i , \delta_i$ are equal and independent of $i$ for 
$1 \leq i \leq d-1$.

We first show $\epsilon_{d-1} = \delta_{1}$.  Setting $(i,j) = (2,1)$ in
(\ref{eq:ttot}), we
find
\begin{eqnarray}
\frac{\theta_{d-2} -\theta_{d}}{\theta_{d-1}-\theta_{d}}=
\frac{\eta_{0}-\eta_{2}}{\eta_{1} -\eta_{2}}.
\label{eq:l1}
\end{eqnarray}
In this equation, the left-hand side (resp. right-hand
side) is equal to $1+\epsilon_{d-1}^{-1}$ (resp. $1+\delta_{1}^{-1}$), so $\epsilon_{d-1}=\delta_1$.
From now on assume $d \geq 3$; otherwise we are done.

We will need the fact that
$\epsilon_{d-1}=\epsilon_{d-2}=\delta_1=\delta_2$.  We already
showed that $\epsilon_{d-1}=\delta_1$.  We will show
$\epsilon_{d-2}=\delta_2$ and $\epsilon_{d-1}=\delta_2$.
We start by showing $\epsilon_{d-2}=\delta_2$.  Setting $(i,j) =
(3,1)$ in (\ref{eq:ttot}) we find
\begin{eqnarray}
\frac{\theta_{d-3} -\theta_{d}}{\theta_{d-1}-\theta_{d}}=
\frac{(\eta_{0}-\eta_{2})}{(\eta_{1} -\eta_{2})}
\frac{(\eta_{0}-\eta_{3})}{(\eta_{1} -\eta_{3})}.
\label{eq:l2}
\end{eqnarray}
Setting $(i,j) = (3,2)$ in (\ref{eq:ttot}) we find
\begin{eqnarray}
\frac{(\theta_{d-3} -\theta_{d})}{(\theta_{d-2}-\theta_{d})}
\frac{(\theta_{d-3} -\theta_{d-1})}{(\theta_{d-2}-\theta_{d-1})}=
\frac{\eta_{0}-\eta_{3}}{\eta_{2} -\eta_{3}}.
\label{eq:l3}
\end{eqnarray}
Multiplying (\ref{eq:l1}) by (\ref{eq:l3}) and dividing
the result by  (\ref{eq:l2}) we find
$$
\frac{\theta_{d-3} -\theta_{d-1}}{\theta_{d-2}-\theta_{d-1}}=
\frac{\eta_{1}-\eta_{3}}{\eta_{2} -\eta_{3}}.
$$
In this equation, the left-hand side (resp. right-hand
side) is equal to $1+\epsilon_{d-2}^{-1}$ (resp. $1+\delta_2^{-1}$), so
$\epsilon_{d-2}=\delta_2$.

We now show $\epsilon_{d-1}=\delta_2$.  Dividing  (\ref{eq:l2}) by 
(\ref{eq:l1}) and subtracting $1$ from the result we find
$$
\frac{\theta_{d-3} -\theta_{d-2}}{\theta_{d-2}-\theta_{d}}=
\frac{\eta_{0}-\eta_{1}}{\eta_{1} -\eta_{3}}.
$$
In this equation, the left-hand side (resp. right-hand
side) is equal to $\epsilon_{d-1}^{-1} \epsilon_{d-2}^{-1}
(1+\epsilon_{d-1}^{-1})^{-1}$ (resp.
$\delta_1^{-1} \delta_2^{-1} (1+\delta_2^{-1})^{-1}$).
From our above comments $\epsilon_{d-1} \epsilon_{d-2}=\delta_1 \delta_2$.
Therefore $1+\epsilon_{d-1}^{-1}=1+\delta_2^{-1}$ so
$\epsilon_{d-1}=\delta_2$.  
We have now shown $\epsilon_{d-1}=\epsilon_{d-2}=\delta_1=\delta_2$.

Applying Lemma \ref{lem:eval} to the eigenvalue sequence $\eta_0,
\eta_1, \ldots, \eta_d$ and since $\delta_1 =\delta_2$ we find
$\delta_i$ is indepentent of $i$ for $1 \leq i \leq d-1$.  
Applying Lemma \ref{lem:eval} to the eigenvalue sequence $\theta_d,
\theta_{d-1}, \ldots, \theta_0$ and since $\epsilon_{d-1}
=\epsilon_{d-2}$ we find
$\epsilon_i$ is independent of $i$ for $1 \leq i \leq d-1$.
Since $\epsilon_{d-1} = \delta_{1}$ we find $\epsilon_i$, $\delta_{i}$
are equal and independent of $i$ for $1 \leq i \leq d-1$. 
\qed

\begin{theorem}
\label{lem:aog}
With reference to Definition \ref{def:s-u}, either $(i)$ or
$(ii)$ below holds.

\noindent 
\begin{enumerate}
\item
Each of the sequences 
$$
\theta_0, \theta_1, \ldots , \theta_d; \qquad \qquad \theta_0^*, \theta_1^*,
\ldots , \theta_d^*; \qquad \qquad \eta_0, \eta_1, \ldots , \eta_d; \qquad
\qquad \eta_0^*, \eta_1^*, \ldots , \eta_d^*
$$
is arithmetic.
\item
There exists $q \in \Bbb K$ such that $q \not= 0$, $q \not= 1$ and
each of the sequences  
$$
\theta_0, \theta_1, \ldots , \theta_d; \qquad \qquad \theta_0^*, \theta_1^*,
\ldots , \theta_d^*; \qquad \qquad \eta_0, \eta_1, \ldots , \eta_d; \qquad
\qquad \eta_0^*, \eta_1^*, \ldots , \eta_d^*
$$
is $q$-classical.
\end{enumerate}
\end{theorem}

{\it Proof.}  Assume $d \geq 2$; otherwise $(i)$ holds trivially.

By Lemma \ref{lem:ca}, $(A^*,A)$ is adjacent to
$(B^*,B)$.  Apply Lemma \ref{lem:vis3} to this pair by replacing
$\theta_i$ with $\theta_i^*$ and $\eta_i$ with $\eta_i^*$.  We find
\begin{eqnarray}
\frac{\theta_{i}^*-\theta_{i+1}^*}{\theta_{i-1}^* -\theta_{i}^*},
\qquad \qquad
\frac{\eta_{i}^* -\eta_{i+1}^*}{\eta_{i-1}^*-\eta_{i}^*}
\label{eq:tinky}
\end{eqnarray}
are equal and independent of $i$ for $1 \leq i \leq d-1$.

By Lemma \ref{lem:ca}, $(A,A^*)$ is adjacent to
$(B^*,B)$.  Apply Lemma \ref{lem:vis3} to this pair by replacing
$\theta_i$ with $\theta_{d-i}$ and $\eta_i$ with $\eta_{d-i}^*$.  We find
\begin{eqnarray}
\frac{\theta_{i}-\theta_{i+1}}{\theta_{i-1} -\theta_{i}}, \qquad \qquad
\frac{\eta_{i}^*-\eta_{i+1}^*}{\eta_{i-1}^* -\eta_{i}^*}
\label{eq:lala}
\end{eqnarray}
are equal and independent of $i$ for $1 \leq i \leq d-1$.

Combining Lemma \ref{lem:vis3} and lines (\ref{eq:tinky}),
(\ref{eq:lala}) we find the scalars
\begin{eqnarray}
\frac{\theta_{i}-\theta_{i+1}}{\theta_{i-1} -\theta_{i}},
 \qquad
\frac{\theta^*_{i}-\theta^*_{i+1}}{\theta^*_{i-1} -\theta^*_{i}},
 \qquad
\frac{\eta_{i}-\eta_{i+1}}{\eta_{i-1} -\eta_{i}},
 \qquad
\frac{\eta_{i}^*-\eta_{i+1}^*}{\eta_{i-1}^* -\eta_{i}^*}
\label{eq:po}
\end{eqnarray}
are equal and independent of $i$ for $1 \leq i \leq d-1$.

Let $q$ denote the common value of (\ref{eq:po}) and observe $q \not= 0$.
If $q = 1$ then $(i)$ holds by Lemma \ref{lem:seq1} and Definition
\ref{def:had}.
If $q \not= 1$ then $(ii)$ holds by Lemma \ref{lem:seq2} and Definition
\ref{def:has}.
\qed

\medskip

\section{Three Mutually Adjacent Leonard Pairs}

We now present our second main result.

\begin{theorem}
Let $(A,A^*)$, $(B,B^*)$, and $(C,C^*)$ denote three mutually adjacent
Leonard pairs on $V$.  Then for each of these pairs, the eigenvalue
sequences and dual eigenvalue sequences are arithmetic.
\end{theorem}

{\it Proof.}  Assume $\dim V \geq 3$; otherwise the result is trivial.

Assume the theorem is false.  Then relabeling the Leonard pairs if
necessary we can assume that the eigenvalue sequences or dual
eigenvalue sequences of $(A,A^*)$ are not arithmetic.  Since the
Leonard pairs $(A,A^*)$ and $(B,B^*)$ are adjacent, we adopt the
notation of Definition \ref{def:s-u}.  Since Theorem
\ref{lem:aog}$(i)$ does not hold for $(A,A^*)$ and $(B,B^*)$,  Theorem
\ref{lem:aog}$(ii)$ must hold.  By Theorem \ref{lem:aog}$(ii)$,
there exists $q \in \Bbb K$ such that $q \not= 0$, $q \not= 1$ and
the eigenvalue sequence for $(A,A^*)$ corresponding to
$[wx]$ and the dual eigenvalue sequence for $(A,A^*)$ corresponding to
$[yz]$ are $q$-classical.

By Lemma \ref{lem:rty}, we find
$\mathcal F  (A,A^*) = \mathcal F  (B,B^*) = \mathcal F  (C,C^*)$;
call this common set $\mathcal F $.  Also by Lemma \ref{lem:rty}, the
principal relations induced on $\mathcal F $ by $(A,A^*)$, $(B,B^*)$,
and $(C,C^*)$ are mutually distinct.
From the construction, the principal relation induced by $(A,A^*)$
partitions $\mathcal F $ into the equivalence classes $\{ w,x \}$ and
$\{ y,z \}$.  Similarly, the principal relation induced by $(B,B^*)$
partitions $\mathcal F $ into the equivalence classes
$\{ w,z \}$ and $\{ x,y \}$.  Therefore the principal relation induced by
$(C,C^*)$ must partition $\mathcal F $ into the equivalence classes
$\{ w,y \}$ and $\{ x,z \}$.  Interchanging $C$ and $C^*$ if
necessary, we can assume that $w,y$ are the $C$-standard flags and
$x,z$ are the $C^*$-standard flags. 

Now apply Theorem \ref{lem:aog} to the adjacent Leonard pairs $(A,A^*)$ and
$(C,C^*)$.  We find Theorem \ref{lem:aog}$(ii)$ holds and there exists
$q' \in \Bbb K$ ($q' \not= 0$, $q' \not= 1$) such that the eigenvalue
sequence for $(A,A^*)$ corresponding to $[wx]$ and the dual eigenvalue
sequence for $(A,A^*)$ corresponding to $[zy]$ is $q'$-classical.
Since the eigenvalue sequence for $(A,A^*)$ corresponding to $[wx]$ is
both $q$-classical and $q'$-classical we have $q=q'$.  Since the
decomposition $[yz]$ is the inversion of the decomposition $[zy]$,
the dual eigenvalue sequence for $(A,A^*)$ corresponding to $[zy]$ is 
both $q^{-1}$-classical and $q'$-classical.  Therefore $q^{-1} =q'$.
Since $q=q'$ and $q^{-1} =q'$, we find $q=1$ or $q=-1$.  By construction
$q \not= 1$.  Observe $q \not= -1$ by Note \ref{n:ote} and since $d
\geq 2$.  We now have a contradiction.
\qed

\section{Example}

In this section we give an example of three mutually adjacent Leonard
pairs.  

Our discussion will start with the Lie algebra $sl_2(\Bbb K)$.
This algebra has a basis $e,f,h$ that satisfies  $[e,f]=h$, $[e,h]=-2e$, and
$[f,h]=2f$, where $[ , ]$ denotes the Lie bracket.  Such a basis is
called a {\it Chevalley basis}.

We recall the irreducible, finite dimensional $sl_2(\Bbb K)$-modules.
For an integer $d \geq 0$, up to isomorphism there exists a unique
irreducible $sl_2(\Bbb K)$-module with dimension $d+1$.  We call this
module $V^d$.  Given a Chevalley basis $e,f,h$ for $sl_2(\Bbb K)$ there
exists a basis for $V^{d}$ with respect to which the matrices
representing $e$, $f$, and $h$ are as follows:
\begin{eqnarray}
e:\left( \begin{array}{cccccc}
0& d &&&&   {\bf 0} \\
&0 & d-1&&&\\
&&0& \cdot            &&\\
&&&\cdot & \cdot& \\
&&&& \cdot & 1 \\
{\bf 0} &&&&& 0 \\
\end{array} \right), 
\qquad \qquad
f:\left( \begin{array}{cccccc}
0& &&&&  {\bf 0} \\
1&0 &&&&\\
&2&0&&&\\
&&\cdot&\cdot && \\
&&&\cdot& \cdot & \\
{\bf 0} &&&&d& 0 \\
\end{array} \right), \label{eq:ef}
\end{eqnarray}
\begin{eqnarray}
h: \diag (d, d-2, \ldots, -d). \label{eq:efh}
\end{eqnarray}

We have a comment about the 2-dimensional irreducible $sl_2(\Bbb K)$-module.

\begin{lemma} 
\label{lem:cheb}
Let $v_0,v_1$ denote a basis for $V^1$.  Then there exists a Chevalley
basis $e,f,h$ for $sl_2(\Bbb K)$ such that
\beast
ev_0=0, \qquad \qquad ev_1=v_0,
\\
fv_0=v_1, \qquad \qquad fv_1=0,
\\
hv_0=v_0, \qquad \qquad hv_1=-v_1.
\eeast
\end{lemma}

{\it Proof.}   Let $e:V^1 \rightarrow V^1$ denote the
linear transformation that sends $v_0$ to $0$ and $v_1$ to $v_0$.  Let
$f:V^1 \rightarrow V^1$ denote the linear transformation that sends
$v_0$ to $v_1$ and $v_1$ to $0$.  Let $h:V^1
\rightarrow V^1$ denote the linear transformation that sends $v_0$ to
$v_0$ and $v_1$ to $-v_1$.  We see that $e,f,h$ have trace $0$
on $V^1$ and therefore can be viewed as elements of $sl_2(\Bbb K)$.
Notice that $e,f,h$ are linearly independent and hence form a basis
for $sl_2(\Bbb K)$.  We check that $[e,f]=h$, $[e,h]=-2e$, and
$[f,h]=2f$.  Therefore $e,f,h$ is a Chevalley basis for $sl_2(\Bbb
K)$.
\qed

\begin{lemma}
\label{lem:ptl}
Let $a,a^* \in sl_2(\Bbb K)$.  Then the following $(i)$--$(iii)$ are
equivalent.
\begin{enumerate}
\item
There exist pairwise linearly independent vectors
$v_0,v_1,w_0,w_1$ in $V^1$ such that
\beast
av_0 = v_0, \qquad \qquad av_1 = -v_1,
\\
a^*w_0 = w_0, \qquad \qquad a^*w_1 = -w_1.
\eeast
\item
$a,a^*$ generate $sl_2(\Bbb K)$, and on $V^1$ we have $\det (a)=\det
(a^*) = -1$.
\item
There exists a Chevalley basis $e,f,h$ for $sl_2(\Bbb K)$ and there
exist $\alpha,\beta,\gamma \in \Bbb K$ such that 
\begin{eqnarray}
a=h, \qquad a^*= \alpha h+\beta e+\gamma f, 
\label{eq:!}
\end{eqnarray} 
and $\beta \gamma = 1-\alpha^2 \not= 0$.
\end{enumerate}
\end{lemma}

{\it Proof.}  $(ii) \Rightarrow (i)$:  Observe that the action of $a$
on $V^1$ has
determinant $-1$ and trace $0$, so the characteristic polynomial of $a$ on
$V^1$ is $x^2-1$.  Therefore the eigenvalues of $a$ on $V^1$ are $1$
and $-1$.  Let $v_0 \in V^1$ denote an eigenvector of $a$
with eigenvalue $1$ and let $v_1 \in V^1$ denote an eigenvector of
$a$ with eigenvalue $-1$.  Note that $v_0,v_1$ are linearly
independent and that $av_0 = v_0$, $av_1 = -v_1$.  Similarly,
there exist linearly independent vectors $w_0,w_1 \in V^1$ such that
$a^*w_0 = w_0$ and $a^* w_1 = -w_1$.  

It remains to show that $v_i,w_j$ are linearly independent for
$i,j \in \{ 0,1 \}$.  Suppose there exist $i,j \in \{ 0,1 \}$ such that
$v_i,w_j$ are linearly dependent.  Then $\spn (v_i) =\spn (w_j)$ is a
proper nonzero subspace of $V^1$ which is invariant under the actions
of $a,a^*$.  This is impossible since $a,a^*$ generate $sl_2(\Bbb K)$
and $V^1$ is irreducible as an $sl_2(\Bbb K)$-module.

$(i) \Rightarrow (iii)$:  By assumption $v_0,v_1$ are linearly
independent so they form a basis for $V^1$.  Let $e,f,h$ denote the
corresponding Chevalley basis for $sl_2(\Bbb K)$ from Lemma
\ref{lem:cheb}.  Comparing
the actions of $a$ and $h$ on the basis $v_0,v_1$ we find $a=h$.
Since $e,f,h$ form a basis for $sl_2( \Bbb K )$ there exist $\alpha,
\beta ,\gamma  \in \Bbb K$ such that $a^*=\alpha h+\beta e+ \gamma f$.
Note that $\beta \not= 0$; otherwise $v_1$ is an eigenvector of $a^*$
and therefore is a scalar multiple of $w_0$ or $w_1$.  Similarly
$\gamma \not= 0$.  Therefore $\beta \gamma \not= 0$.  We show $ \beta
\gamma = 1 - \alpha^2$.  Observe $w_0,w_1$ is a basis for $V^1$.  From the
action of $a^*$ on this basis we find the determinant of $a^*$ on
$V^1$ is $-1$.
Using $a^*=\alpha h+\beta e+ \gamma f$ and the data in Lemma
\ref{lem:cheb} we find the determinant of $a^*$ on $V^1$ is $-
\alpha^2 - \beta \gamma$.  Therefore $-1 = - \alpha^2 - \beta \gamma$
so $ \beta \gamma = 1 - \alpha^2$.

$(iii) \Rightarrow (ii)$:  Since $e,f,h$ is a Chevalley basis for
$sl_2(\Bbb K)$, there exists a basis $v_0,v_1$ for $V^1$ such that
$ev_1=v_0$, $ev_0=0$, $fv_0=v_1$, $fv_1=0$, $hv_0=v_0$, $hv_1=-v_1$.
With respect to the basis $v_0,v_1$ the matrices representing $a,a^*$ are
\begin{eqnarray}
a:\left( \begin{array}{cc}
1 & 0 \\
0 & -1
\end{array} \right), 
\qquad \qquad
a^*:
\left( \begin{array}{cc}
\alpha & \beta \\
\gamma & -\alpha
\end{array} \right).
\label{eq:bnm}
\end{eqnarray}
We show $a,a^*$ generate $sl_2(\Bbb K)$.  From (\ref{eq:bnm}) we find
$[a,a^*] = 2\beta e-2\gamma f$.  Comparing this with (\ref{eq:!}) we
find $a,a^*, [a,a^*]$ are linearly independent and therefore span
$sl_2(\Bbb K)$.  This shows $a,a^*$ generate $sl_2(\Bbb K)$.  From
(\ref{eq:bnm}) we find that on $V^1$ we have $\det (a)=\det (a^*)=-1$.
\qed

\begin{lemma}
\label{lem:actlp}
Let $a, a^* \in sl_2(\Bbb K)$ satisfy the equivalent conditions
$(i)$--$(iii)$ in Lemma \ref{lem:ptl}.  Let $d$ denote a nonnegative
integer.  Then $a,a^*$ act on $V^{d}$ as a Leonard pair.  The sequence
$d-2i$ ($0 \leq i \leq d$) is both an eigenvalue sequence and a dual
eigenvalue sequence for this pair.
\end{lemma}

{\it Proof.} 
By Lemma \ref{lem:ptl}$(iii)$ there exists a Chevalley basis $e,f,h$
for $sl_2(\Bbb K)$ and there exist $\alpha,\beta,\gamma \in \Bbb K$
with $\beta,\gamma$ nonzero such that
$a=h$ and $a^*=\alpha h+\beta e+\gamma f$.  Consider the basis for
$V^d$ with respect to which the matrices representating $e,f,h$
are given in (\ref{eq:ef}) and (\ref{eq:efh}).  With respect to this
basis, the matrices representating $a$ and $a^*$ are: 
\begin{eqnarray}
a:\diag (d,d-2, \ldots, -d), 
\end{eqnarray}
$$
a^*:\left( \begin{array}{ccccccc}
d \alpha &d \beta &&&&&  {\bf 0} \\
\gamma&(d-2) \alpha & (d-1)\beta&&&&\\
&2\gamma&(d-4) \alpha& (d-2)\beta&&&\\
&&3\gamma&\cdot & \cdot&& \\
&&&\cdot& \cdot & \cdot&\\
&&&&\cdot &\cdot&\beta\\
{\bf 0} &&&&&d \gamma & -d \alpha  \\
\end{array} \right) .
$$
Since $\beta, \gamma$ are nonzero, the matrix representing $a^*$ is
irreducible tridiagonal.

Replacing $(a,a^*)$ by $(a^*,a)$ in the argument so far, we find there
exists a basis for $V^d$ with respect to which $a^*$ acts as
$\diag (d,d-2, \ldots, -d)$ and
$a$ is irreducible tridiagonal.  The result follows.
\qed

\begin{definition}
\label{def:3lps}
\rm
Let $v_0,v_1,w_0,w_1$ denote pairwise linearly independent vectors in
$V^1$.  We let $a,a^*,b,b^*,c,c^*$ denote the elements of $sl_2(\Bbb
K)$ that satisfy the following:
\beast
av_0 = v_0, \qquad \qquad av_1 = -v_1,
\\
a^*w_0 = w_0, \qquad \qquad a^*w_1 = -w_1,
\\
bv_0 = v_0, \qquad \qquad bw_0 = -w_0,
\\
b^*w_1 = w_1, \qquad \qquad b^*v_1 = -v_1,
\\
cv_0 = v_0, \qquad \qquad cw_1 = -w_1,
\\
c^*w_0 = w_0, \qquad \qquad c^*v_1 = -v_1.
\eeast
\end{definition}

\medskip

\begin{theorem}
\label{th:3lps}
Let $d$ denote a nonnegative integer.  With reference to
Definition \ref{def:3lps}, the pairs $(a,a^*)$, $(b, b^*)$,
$(c, c^*)$ act on $V^{d}$ as mutually adjacent Leonard pairs.
\end{theorem}

{\it Proof.}  Abbreviate $V = V^d$.  By Lemma \ref{lem:actlp} each of
the pairs $(a,a^*)$, $(b, b^*)$, $(c, c^*)$ acts on $V$ as a Leonard
pair.  We show that these Leonard pairs are mutually adjacent.  We start by
showing the first two are adjacent.

By assumption $v_0,v_1$ are linear
independent so they form a basis for $V^1$.  Let $e,f,h$ denote the
corresponding Chevalley basis from Lemma \ref{lem:cheb}.  Observe that
$a=h$.  We write each of $a^*,b,b^*$ as a linear combination of
$e,f,h$ and consider the corresponding coefficients.
For $a^*$, the coefficient of $e$ (resp. $f$) is nonzero; otherwise
$v_1$ (resp. $v_0$) is an eigenvector of $a^*$ and therefore is a
scalar multiple of $w_0$ or $w_1$.  For $b$, the coefficient of
$e$ is nonzero since $v_1$ is not an eigenvector of $b$ and the
coefficient of $f$ is zero since $v_0$ is an eigenvector of $b$.  For
$b^*$, the coefficient of $e$ is zero since $v_1$ is an eigenvector
of $b^*$ and the coefficient of $f$ is nonzero since $v_0$ is not an
eigenvector of $b^*$.  Now consider the basis for $V$ with respect
to which the matrices representing $e,f,h$ are given by (\ref{eq:ef})
and (\ref{eq:efh}).  By our above comments, with respect to this basis
$a$ is diagonal, $a^*$ is irreducible tridiagonal, $b$ is upper
bidiagonal, and $b^*$ is lower bidiagonal.  
Let $V_0,V_1, \ldots, V_d$ denote the decomposition of $V$ induced by
this basis.  Then $V_0,V_1, \ldots, V_d$ is $a$-standard and
$UL$-split for $(b,b^*)$.  Observe the decomposition $V_d,V_{d-1},
\ldots, V_0$ is $a$-standard and $LU$-split for $(b,b^*)$.  Given a
decomposition of $V$ that is $a$-standard, this decomposition is
either $V_0,V_1, \ldots, V_d$ or $V_d,V_{d-1},
\ldots, V_0$.  In either case, this decomposition is split for
$(b,b^*)$.  We have now shown that each $a$-standard decomposition of
$V$ is split for $(b,b^*)$.  In a similar fashion, we find that each
$a^*$-standard decomposition of $V$ is split for $(b,b^*)$.
Now each decomposition of $V$ that is standard for $(a,a^*)$ is split
for $(b,b^*)$.  Therefore $(a,a^*)$ and $(b,b^*)$ satisfy Lemma
\ref{th:doalp}$(ii)$.  Applying Definition \ref{def:ad1}, we find the
actions of $(a,a^*)$ and $(b,b^*)$ on $V$ are adjacent Leonard pairs.

The proofs that $(a,a^*)$ is adjacent to $(c,c^*)$ and $(b,b^*)$ is
adjacent to $(c,c^*)$ are similar and are left to the reader.
\qed

\section{Leonard Pairs With Arithmetic Eigenvalue And Dual Eigenvalue
  Sequences} 

In this section we will show that if $(A,A^*)$ is a Leonard pair on
$V$ with arithmetic eigenvalue and dual eigenvalue sequences, then there
exist Leonard pairs $(B,B^*)$ and $(C,C^*)$ on $V$ such that $(A,A^*)$,
$(B,B^*)$ and $(C,C^*)$ are mutually adjacent.  Before we present our
result, we will first discuss the notion of an {\it affine
  transformation} of a Leonard pair. 

Let $(A,A^*)$ denote
a Leonard pair on $V$.   Let $\alpha, \beta, \alpha^*, \beta^*$ denote
scalars in $\Bbb K$ such that $\alpha \not= 0$, $\alpha^* \not= 0$.
One easily verifies that
\begin{eqnarray}
(\alpha A+\beta I,\alpha^* A^*+\beta^* I)
\label{eq:aftr}
\end{eqnarray}
is a Leonard pair on $V$.  We call (\ref{eq:aftr}) an {\it affine
transformation} of $(A,A^*)$.  Observe that $(A,A^*)$ and
(\ref{eq:aftr}) have the same eigenspace and dual eigenspace
decompositions.  Therefore a Leonard pair on $V$ is adjacent to
$(A,A^*)$ if and only if it is adjacent to (\ref{eq:aftr}).  Let
$\theta_i$ (resp. $\theta_i^*$) ($0 \leq i \leq d$) denote an
eigenvalue sequence (resp. dual eigenvalue sequence) of $(A,A^*)$.  Then 
$\alpha \theta_i + \beta$ (resp. $\alpha^* \theta_i^* + \beta^*$)
($0 \leq i \leq d$) is an eigenvalue sequence (resp. dual eigenvalue
sequence) of (\ref{eq:aftr}).

We now present our third main result.

\begin{theorem}  
\label{th:last}
Let $(A, A^*)$ denote a Leonard pair on $V$
with an arithmetic eigenvalue sequence and an arithmetic dual
eigenvalue sequence.   Then there exist Leonard
pairs $(B,B^*)$ and $(C,C^*)$ on $V$ such that $(A,A^*)$, $(B,B^*)$ and
$(C,C^*)$ are mutually adjacent.
\end{theorem}

{\it Proof.}
Let $d= \dim V-1$.  Applying an affine transformation to $(A,A^*)$ if
necessary, we can assume that $d-2i$ ($0 \leq
i \leq d$) is an eigenvalue sequence and a dual eigenvalue sequence of
$(A,A^*)$.  We will first show that there exist $a,a^* \in sl_2(\Bbb
K)$ satisfying the equivalent conditions $(i)$--$(iii)$ of Lemma
\ref{lem:ptl} such that the action of $a, a^*$ on
$V^{d}$ is a Leonard pair isomorphic to $(A,A^*)$.

By \cite[Example 5.13]{ter7}, we see that $(A,A^*)$ is of Krawtchouk
type; now by \cite[Theorem 9.1, Theorem 9.3]{ter7}, there exists a
scalar $p \in \Bbb K$ ($p \not= 0, p \not= 1$) and there exists a
basis for $V$ with respect to which the matrices representing $A$ and
$A^*$ are  
\begin{eqnarray}
A: \diag (d,d-2, \ldots, -d),
\label{eq:matA}
\end{eqnarray}
\begin{eqnarray}
A^*:\left( \begin{array}{cccccc}
\alpha_0 &\beta_0 &&&&  {\bf 0} \\
\gamma_1& \alpha_1& \beta_1&&&\\
& \gamma_2&\cdot & \cdot&& \\
&&\cdot& \cdot & \cdot&\\
&&&\cdot &\cdot& \beta_{d-1}\\
{\bf 0} &&&& \gamma_d & \alpha_d \\
\end{array} \right),
\label{eq:woohoo}
\end{eqnarray}
where $\alpha_i = (1-2p)(d-2i)$, $\beta_i =
2p(d-i)$, and $\gamma_i = 2(1-p)i$.

Let $e,f,h$ denote a Chevalley basis for $sl_2(\Bbb K)$.
Define $a=h$ and $a^*=(1-2p)h +2pe +2(1-p)f$.  Notice that $a,a^*$
satisfy Lemma \ref{lem:ptl}$(iii)$, so by Lemma \ref{lem:actlp},
$(a,a^*)$ acts as a Leonard pair on $V^{d}$.
Consider the basis for $V^{d}$ with respect to which the matrices
representing
$e,f,h$ are given in (\ref{eq:ef}) and (\ref{eq:efh}).  With respect
to this basis the matrices representing $a$ and $a^*$ are
given in (\ref{eq:matA}) and (\ref{eq:woohoo}) respectively.  It is now
apparent that the Leonard pair $(A,A^*)$ is isomorphic to the Leonard
pair $(a,a^*)$ on $V^d$.
Because of this, it suffices to show the Leonard pair $(a,a^*)$ on
$V^d$ is part of three mutually adjacent Leonard pairs.  From Lemma
\ref{lem:ptl}$(i)$, there exist pairwise linearly independent vectors
$v_0,v_1, w_0,w_1$ in $V^1$ such that
\beast
a v_0 = v_0, \qquad \qquad a v_1 = -v_1,
\\
a^* w_0 = w_0, \qquad \qquad a^* w_1 = -w_1.
\eeast
Define $b,b^*,c,c^* \in sl_2(\Bbb K)$ as in Definition
\ref{def:3lps}.  By Theorem \ref{th:3lps}, $(a,a^*)$, $(b,b^*)$,
$(c,c^*)$ act on $V^d$ as mutually adjacent Leonard pairs.  The result
follows.
\qed

\bigskip
\bigskip
\bigskip

\noindent
{\bf Acknowledgements:}  This paper was written while the author was
a graduate student at the University of Wisconsin-Madison.  The author
would like to thank his advisor Paul Terwilliger for his many valuable
ideas and suggestions.  The author would also like to thank Brian
Curtin, Eric Egge and Arlene Pascasio for giving this manuscript a
careful reading and offering many valuable suggestions.

\noindent Brian Hartwig  \hfil \break
Department of Mathematics  \hfil \break
University of Wisconsin  \hfil \break
480 Lincoln Drive \hfil \break
Madison, Wisconsin, 53706, USA \hfil \break
email: hartwig@math.wisc.edu \hfil \break

\end{document}